\documentclass[final]{IEEEtran}

\newcommand{\Real}[0]{\mathbb R}

\newcommand{\Natural}[0]{\mathbb N}

\renewcommand{\exp}[1]{\mathit{e}^{#1}}

\newcommand{\conj}[1]{{#1}^*}

\newcommand{\T}[0]{T}
\newcommand{\mv}[1]{\mathbf{#1}}

\newcommand{\td}[1]{\dot{#1}}
\newcommand{\Laplace}{\Delta}

\usepackage{amsmath}

\usepackage{graphicx}
\usepackage{amssymb}
\usepackage{ae}
\usepackage{amsfonts}
\usepackage{newlfont}
\usepackage{epstopdf}

\title{Probabilistic tractography, Path Integrals and the Fokker Planck equation}
\author{Marco Reisert \\
Dept. of Radiology, Medical Physics, University Hospital Freiburg, Germany}

\begin{document}

\maketitle

\begin{abstract}
Probabilistic tractography based on diffusion weighted MRI has become a powerful approach for quantifying structural 
brain connectivities. In several works the similarity of probabilistic tractography and path integrals was already pointed out. This work
investigates this connection more closely. For the so called Wiener process, a Gaussian random walker,
the equivalence is worked out. We identify the source of the asymmetry of usual random walkers approaches 
and show that there is a proper symmetrization, which leads to a new symmetric connectivity measure. To compute this measure we will
use the Fokker-Planck equation, which is an equivalent representation of a Wiener process in terms of a partial differential 
equation. In experiments we show that the proposed approach leads a symmetric and robust connectivity measure. 
\end{abstract}

\section{Introduction}
Diffusion MRI has become a very important tool for understanding the living brain tissue (\cite{Jones2010_diff_book}).
It can reveal both macro- and microstructural features 
of the neuronal network of the human brain. Tractography tries to characterize the structural connectome to understand the details of 
the interregional relationships of the human brain. Tractography algorithms may be divided in deterministic, streamline-based 
methods (\cite{Mori99,Basser2000}), probabilistic approaches (\cite{Hagmann2003,Parker2003}) and global approaches (\cite{Mangin2002,reisert2010}).
In probabilistic tractography one basically draws samples from a distribution over paths and 
computes some statistics over these paths, for example, recording their endpoints. 
In some work the similarity to the notion of path integrals appearing 
in quantum mechanics \cite{feynman2012quantum} and statistical physics \cite{kleinert2009path} was already pointed out \cite{tuch2000path,bjornemo2002regularized,friman2006bayesian,Behrens2007}. Rigorous mathematical 
investigations show that the basic stochastic process behind path integrals is 
a so called Wiener process, a continuous Gaussian random walker \cite{vanKampen1992stochastic}. 
In this work we will recap the foundations of the theory of Wiener processes and path integrals. 
Based on this we build a path integral that leads to a symmetric brain connectivity measure. We will see that the 
source of asymmetry of conventional walker principles is due to particle conservation. By dropping particle conservation and starting from the path integral perspective we will find a symmetric brain connectivity measure. Besides the walker perspective and path integrals, there is a third, equivalent 
approach, which describes expectation values of Wiener walkers by partial differential equations (PDE). This equivalence will be used to give an algorithm 
for the computation of the connectivity measure. It will be based on solving a large linear system, describing a mixed diffusion/convection process. 
We also introduce a novel discretization scheme to avoid the heavy directional dependency, which usually appears for discretized convection operators. 
PDEs in the context of diffusion MRI have been previously used for regularization \cite{reisert2011}, 
connectivity estimation \cite{tournier2003diffusion,hageman2009diffusion,zhang2013logical} and fiber density estimation \cite{miccai12}.
Our proposed PDE contains a convection operator, but in the joint position/orientation space, in \cite{batchelor2002fibre} also convection was
used for tractography, however just in position space. 
In \cite{zhang2013logical} probabilistic tractography was put into a logical framework 
and and an algorithm was derived which is based on solving a large linear system. This approach also works just in position space such that 
crossings cannot handled adequately as long as only local neighborhoods are considered. The system matrix describes a anisotropic diffusion process.
%
%
%
%
%
%
%
%
%
%
%
%
%
%
%

\section{Method}
Apart from a few examples most methods to estimate brain connectivity are based on the walker principle. 
Fiber tracts are initiated from certain seed points and are iteratively built by following locally defined directions. 
While the deterministic tracking approaches are more used 
for illustrative purposes, the probabilistic ones are more related to quantitative connectivity analysis. The mathematical
principle behind both approaches is an integration of streamlines along the underlying data. In probabilistic tracking
the process can be seen as a Markov process. The iteration process is simple: if we call $\mv s(t)$ the current state of the tracker at step $t$,
then the next state $\mv s(t+1)$ is drawn from some transition probability density
$W(\mv s(t+1) | \mv s(t))$, which may depend in some way on the DW-measurement. This process is basically a random walk in the state space
of the tracker. If $W$ is Gaussian, it is possible to formulate the limit for very small time steps, which results in a Wiener process, which
we will now concentrate on. 

\subsection{Wiener Process}
In physics and mathematics there exists a huge collection of concepts for the analysis and characterization of Wiener random walkers. 
Basically, there are three perspectives which are mostly equivalent. 

The first one, which is closest to the above described Markov process 
is the concept of stochastic calculus and stochastic differential equations (SDE) \cite{vanKampen1992stochastic}.
In physics, SDEs are usually written as Langevin equations, which do not have the strong mathematical footing of a SDE,
but are simpler to read and more similar to ordinary calculus.
To give an example (actually this example already covers everything we need later) 
we want to consider a simple diffusion process with some additional drift.
Let the state of the 'tracker' be $\mv s \in \Real^3$ and $\mv v$ a vector field causing the drift. Then, the corresponding
Langevin equation is 
 \begin{equation}
 \td{\mv s} (t) = \mv v(\mv s(t) ) +  \eta(t) \label{eq:langevin}
 \end{equation}
where $\eta$ stands for mean free white noise with variance $\sigma^2$ per unit time. It is uncorrelated in time, 
meaning $\langle\eta(t)\eta(t')\rangle = \sigma^2 \delta(t-t')$.
The dot in $\td{\mv s}$ means time differentiation. 
 An approximative numerical integration gives rise to the 
following propagation scheme:
\begin{equation}
 \mv s(t+\Delta t) = \mv s(t) + \mv v(\mv s(t)) \Delta t +  \mv u(t) \sqrt{\Delta t}.  \label{eq:langdisc}
\end{equation}
The first term $\mv v(\mv s(t)) \Delta t$ looks like ordinary Euler integration. The second, stochastic 
term $\mv u(t)$ is drawn independently for each time step from $\mathcal{N}(0,\sigma^2)$. The factor $\sqrt{\Delta t}$
expresses the fact that
the process $\eta(t)$ is not differentiable: variances add up and not standard deviations. Note the difference between 
the continuous Wiener process $\eta(t)$ and the discrete process $\mv u(t)$ which is just defined at discrete time points $k\Delta t$. 
Both are related by $\int_{t-\Delta t}^t \eta(t) dt = \mv u(t) \sqrt{\Delta t}$.
The transition probability $W$ as described in the beginning for the discrete process described in \eqref{eq:langdisc}
is a Gaussian $W(\cdot|\mv s(t)) = \mathcal{N}(\mv s(t) + \Delta t\  \mv v(\mv s(t)),\Delta t \sigma^2)$.

The second perspective describes expectation values 
of an ensemble $\mathcal{S}$ of random walkers described by \eqref{eq:langevin}. Suppose we have generated such 
an ensemble $\mathcal S$ of  walkers all starting at $\mv s(0) = \mv x_0$ and we want to know the distribution 
of the states at some time $T$, i.e. 
\begin{equation}
p_T(\mv x|\mv x_0) = \sum_{\{\mv s\} \in \mathcal S} \delta(\mv s(T) - \mv x) \label{eq:pdens}
\end{equation}
Note, that the sum ranges over complete random walks or paths. To make this formally more transparent
we always write curly brackets whenever we refer to a path (a chain of states of a walker) as a whole and not 
just a particular state at a specific time. 
In fact, the distribution described in \eqref{eq:pdens} is a solution of a partial differential equation, which is called the Fokker-Planck equation (FPE) and is the master equation 
of the proposed continuous stochastic process. For this example, the FPE takes the form
\begin{equation}
 \td{p_t} = -\nabla \cdot (\mv v p_t) + \frac{\sigma^2}{2} \Laplace p_t  = \mathcal{H} p_t \label{eq:fpe}
\end{equation}
where $\nabla$ is the usual gradient operator in $\Real^3$ and $\Delta$ the Laplacian. If we integrate this 
equation with respect to the initial condition $p_0(\mv x|\mv x_0) = \delta(\mv x - \mv x_0)$ we just resemble the 
distribution given in \eqref{eq:pdens}. The function $p_T(\mv x|\mv x_0)$ is also known as the propagator or Green's
function of the corresponding stochastic process. It can formally be written as $p_t = \exp(\mathcal{H}t)$. Note that $p_T(\mv x|\mv x_0)$ does not necessarily need to be normalized like a probability, e.g. walkers can die at the boundaries, which manifests in the boundary conditions of the corresponding FPE.
The proposed algorithm will be based on discretized solutions of the steady state solutions of a symmetrized FPE.

The third perspective is the path integral concept. From the viewpoint of a brain connectivity, the theory of path integrals is
probably the most appealing one, however, they are probably the mathematical least understood concept 
and do not give a constructive way for designing an algorithm. However, it is essential for the understanding of 
what we are actually doing. The idea is to compute
$p_T(\mv x|\mv x_0)$ by summing over \emph{all} path starting at $\mv x_0$ and ending at $\mv x$  weighted by its probability. 
For a rough derivation and motivation recall the discrete time representation. To compute the total probability that a particular path is drawn we take the product along the path
\[
P(\mv s(N\Delta t),\hdots,\mv s(\Delta t)|\mv s(0) = \mv x_0) = \prod_{i=1}^N W(\mv s(i \Delta t)|\mv s((i-1) \Delta t),
\]
but note that $W$ is a probability \emph{density}. So, it is not totally correct to just multiply them, 
the problem will be discussed below. By taking the logarithm of the above product it can 
be converted into a sum and in the continuous time limit $N \rightarrow \infty$ with $N \Delta t =T$ we have
the integral form
\[
  -\log(P(\{\mv s\}|\mv s(0) = \mv x_0))  = -\int_0^T \log(W(\mv s(t),\td{\mv s}(t)) dt 
\]
The functional $L(\mv s(t),\td{\mv s}(t)) = -\log(W(\mv s(t),\td{\mv s}(t)))$ is sometimes called the Onsager-Machlup functional (\cite{haken1976generalized,durr1978onsager,ito1981characterization,andersson1999finite}) or,
in the style of mechanics, the Lagrangian. It describes the cost of a path. 
Let us further follow the naive way and use \eqref{eq:langdisc} and the 
corresponding $W$ to compute $L$. By disregarding the normalization constant (which actually diverges in the limit),
we can show that 
\begin{equation}
 L_\text{sym}(\mv s,\td{\mv s}) = \frac{1}{2\sigma^2}(\td{\mv s} - \mv v(\mv s))^2 \label{eq:sym}
\end{equation}
which seems to be plausible but is wrong. The correct answer is
\begin{equation}
 L(\mv s,\td{\mv s}) = \frac{1}{2\sigma^2}(\td{\mv s} - \mv v(\mv s))^2 + \frac{1}{2} (\nabla \cdot \mv v)(\mv s) \label{eq:nosym}
\end{equation}
The error is coming from the faulty assumption that we can simply multiply the conditional transition
densities $W$ to get the total probability, which leads to the diverging normalization constant and the 
missing term $\nabla\cdot \mv v$. Actually we have to consider some small volumes instead of 
infinitesimal points. We have to consider the portion of walkers starting in a volume around $\mv s((i-1)\Delta t)$ 
and arriving in a volume around $\mv s(i\Delta t)$. This leads to an additional factor $(1+\frac{\Delta t\ \nabla \cdot \mv v}{2})$
which is the infinitesimal volume change during time $\Delta t$ caused by the drift $\mv v$. In more rigorous 
derivations of this result the factor is justified by the Jacobian of a change of variables (\cite{chaichian2001path,kleinert2009path}),
or by the considering the most probable tube around the path (\cite{durr1978onsager}). 
 
 Let us consider the final equation stating that the probability for a walker starting at $\mv s(0) = \mv x_0$ 
and arriving at $\mv s(T) = \mv x$ is the sum over all paths connecting this two points weighted by their probabilities
\begin{equation}
 p_T(\mv x|\mv x_0) = \sum_{\{\mv s\} \in \mathcal{P}_T^{(\mv x|\mv x_0)}} \exp \left(-\int_0^T L(\mv s,\td{ \mv s}) dt\right)  \label{eq:pi}
\end{equation}
where $\mathcal{P}_T^{(\mv x|\mv x_0)}$ denotes the set of all paths of length $T$ starting in $\mv x_0$ and terminating at $\mv x$.
In fact, the \emph{sum} over all path is an integration in a functional space, for more details see e.g. \cite{feynman2012quantum,kleinert2009path,vanKampen1992stochastic}. For convenience we omit to write here any normalization constant and assume that the measure is appropriately normalized.

\subsection{Reversibility and Symmetry}
In fact, the additional term $(\nabla \cdot \mv v)/2$ in \eqref{eq:nosym} is quite important for us, because it is responsible 
for the fact that the described random walk is irreversible (\cite{garrido1979stochastic}).
In the context of brain connectivity
this is directly related to the symmetry of the connectivity measure, so we will detail this out now.
If we call $L'$ the Lagrangian with negative velocities $\mv v \mapsto -\mv v$ and $\mv s'(t) = \mv s(T-t)$
the path traversed in opposite direction. Then, it is easy to see that 
$\int L'(\mv s',\td{\mv s}') dt \neq \int L(\mv s,\dot{\mv s}) dt$ with the  $L$ obtained in \eqref{eq:nosym}.
This also leads to $p_T(\mv x|\mv x_0) \neq p'_T(\mv x_0|\mv x)$. This asymmetry is basically caused by 
drift together with particle conservation. Considering the corresponding FPE in \eqref{eq:fpe} makes it more clear.
While the diffusion term $\Laplace$ is symmetric (selfadjoint, see Appendix \ref{app1} for definition), the convection/drift is not antisymmetric, i.e.
the operator $\mathcal{H}'$ corresponding to the flipped Lagrangian $L'$ is not equal to adjoint
operator $\mathcal{H}^+$ which has to be case to fulfill $p_T(\mv x|\mv x_0) = p'_T(\mv x_0|\mv x)$.
However, for $L_\text{sym}$ everything is different. In fact, we
have 
\[
 \int_0^T L_\text{sym}'(\mv s',\td{\mv s}') dt= \int_0^T L_\text{sym}(\mv s,\dot{\mv s}) dt
\]
which should also lead to a symmetric propagator $p_T(\mv x|\mv x_0) = p'_T(\mv x_0|\mv x)$.
But how to compute \eqref{eq:pi} for $L_\text{sym}$ in practice, what is the corresponding FPE?  
The term $\nabla \cdot \mv v$ acts like an additional potential energy. In fact, the Feynman-Kac (FK)
formula \cite{kac1949distributions} tells us that such a potential can be directly integrated into the FPE: suppose
we have a given Lagrangian $L$ together with its path integral and a corresponding FPE with
operator $\mathcal{H}$. For a potential field $V(\mv s)$
we define a new $L'(\mv s ,\td{\mv s}) = L(\mv s ,\td{\mv s})-V(\mv s)$, then the FK formula tells us that 
the corresponding FPE operator is just $\mathcal H' = \mathcal H + V$. Applying this our problem with 
$V =  (\nabla \cdot \mv v)/2$ gives
 \begin{eqnarray*}
  \mathcal{H}_\text{sym}p &=& \mathcal{H}p + \frac{1}{2}(\nabla \cdot \mv v)p \\
 &=&   \frac{\sigma^2}{2} \Laplace p + \frac{1}{2}(\nabla \cdot \mv v)p - \nabla \cdot (\mv v p) \\
 &=&   \frac{\sigma^2}{2} \Laplace p - \frac{1}{2}(\mv v \cdot \nabla p + \nabla \cdot (\mv v p) ) 
 \end{eqnarray*}
From the last line we can also see the symmetry property of $\mathcal{H}_\text{sym}$. The 
 convection part $\nabla \cdot (\mv v p) + \mv v \cdot \nabla p$ is antisymmetric in the sense of the operator adjoint and 
hence the operator $\mathcal{H}'_\text{sym}$, which corresponds to flipped velocities,
is just the adjoint of $\mathcal{H}_\text{sym}$. 
One can also see that the particle conservation is lost: we cannot write $\mathcal{H}_\text{sym}$ as the a divergence of something. 
Note the similarity of $\mathcal{H}_\text{sym}$ with the quantum mechanical Hamiltonian of a particle in a vector potential, which is up to 
some constants $(i \nabla - \mv v)^2 = -\Laplace - i(\nabla \mv v - \mv v \nabla) + \mv v^2$. In quantum mechanics the Hamiltonian is
indeed self-adjoint due to the additional complex unit in front of the anti-symmetric part.

\subsection{Steady States and Path Trails}
Up to now we have considered only paths with some specific length $T$. However, we are interested
in \emph{all} path without any length restriction. That is, we sum up over all paths connecting $\mv x_0$
and $\mv x$ of arbitrary length to get a density $p(\mv x | \mv x_0)$ which is
independent of $T$, i.e.
\begin{eqnarray*}
 p(\mv x| \mv x_0) &=& \int_0^\infty dT\  p_T(\mv x| \mv x_0) \\ 
 &=&  \sum_{\{\mv s\} \in \mathcal{P}^{(\mv x|\mv x_0)}} e^{-\int_0^T L(\mv s,\td{ \mv s}) dt}  \\
\end{eqnarray*}
where $\mathcal{P}^{(\mv x|\mv x_0)}$ is the set of paths of any length connecting $\mv x_0$ with $\mv x$. With the assumption $\lim_{T\rightarrow \infty } p_T(\mv x|\mv x_0) = 0$, i.e. all walkers eventually
die, the function $p(\mv x | \mv x_0)$ is the steady state solution of the
 corresponding FPE, i.e. $p(\mv x|\mv x_0)$ is the solution of the equation 
\[
-\mathcal{H} p(\mv x | \mv x_0) = \delta (\mv x-\mv x_0).
\]
This is the basic type of equation we will solve to estimate brain connectivities. 
In the following we call $p(\mv x|\mv x_0)$ the connectivity amplitude. From ordinary probabilistic tractography we know
the so called length bias. With increasing euclidean distance, connectivity values decrease dramatically, usually exponentially.
To account for this bias a linear reweighing was proposed in \cite{Behrens2007}, that is,
\begin{equation}
p_\text{lin} (\mv x| \mv x_0) = \int_0^\infty dT\  T\  p_T(\mv x| \mv x_0). \label{eq:linrew}
\end{equation}
We will see later how $p_\text{lin}$ can be computed in practice. 
Due to the exponential decrease one could also tend to use exponential corrections, i.e. use 
\begin{equation}
p_\kappa(\mv x| \mv x_0) = \int_0^\infty dT\ e^{\kappa T}  p_T(\mv x| \mv x_0)      \label{eq:exprew}                                                 
\end{equation}
with positive $\kappa\in\Real$. In fact, this results in a spectral shift of the operator $\mathcal H \mapsto \mathcal H + \kappa$.
Of course, the choice of $\kappa$ is difficult. Large $\kappa$ will cause bad condition numbers for $\mathcal{H}$ and, in the extreme
case, lead to a diverging $p_\kappa(\mv x| \mv x_0)$.

Until now, we were interested in a measure how strong two points $\mv x_1$ and $\mv x_0$ are connected, but one 
might also be interested in \emph{how} the actual connection looks like. Instead of just summing over all paths connecting 
two points, one can additionally image the path itself by counting how often the path visited a specific position $\mv x$. 
We call this image the trail image of a path $\mv s$, and it is defined by 
\[
 \tau_{\mv s}(\mv x) =  \int_0^T \delta(\mv x - \mv s(t)) dt 
\]
Now, we can write the mean trail image of a path connecting two point $\mv x_1$ and $\mv x_2$ by the expectation value
\[
\bar{\tau}_{(\mv x_1| \mv x_0)} (\mv x) =  \sum_{\{\mv s\} \in \mathcal{P}^{(\mv x|\mv x_0)}}\tau_{\mv s}(\mv x) e^{-\int_0^T L(\mv s,\td{ \mv s}) dt}.
\]
In order to compute $\bar{\tau}$ we have to use the so called Einstein-Smoluchowski-Kolmogorov-Chapman (ESKC) relation, or 
simply semigroup property, which tells that the path-integral, or the corresponding propagator can always be split like
\[
 p_T(\mv x_1|\mv x_0) = \int d\mv x'\ p_{T-t'}(\mv x_1|\mv x') p_{t'}(\mv x'|\mv x_0)
\]
for any given intermediate point $0<t'<T$. This means we can slice the propagator at any given length $t'$ and consider all 
possible intermediate positions $\mv s(t') = \mv x'$ the path has at $t'$ and integrate over them by accounting for the probability
of the path to get from $\mv x_0$ to $\mv x'$  and to reach the target $\mv x_1$ starting from $\mv x'$. 
To compute $\bar{\tau}$ we first consider the point density of the path at some specific length $t'$. That is, we count 
how often a path traverses the point $\mv x$ at length $t'$ by
\begin{eqnarray*}
 \rho_{T,t'}(\mv x) &=& \int d\mv x'\ \delta(\mv x -\mv x') p_{T-t'}(\mv x_1|\mv x') p_{t'}(\mv x'|\mv x_0) \\ 
 &=&  p_{T-t'}(\mv x_1|\mv x) p_{t'}(\mv x|\mv x_0)
\end{eqnarray*}
If we integrate now $\rho_{t',T}$ over all intermediate times $0<t'<T$ and, again, integrate this over all possible $0<T<\infty$ 
we get the mean trail image
\[
\bar{\tau}_{(\mv x_1| \mv x_2)} (\mv x) = \int_0^\infty dT \int_0^T dt' \rho_{T,t'}(\mv x)
\]
Plugging this altogether gives the simple final result
\[
 \bar{\tau}_{(\mv x_1| \mv x_0)} (\mv x) = p(\mv x_1 | \mv x) p(\mv x | \mv x_0)
\]
and, as we already know how to compute $p$ we are ready. Further note that
\begin{eqnarray*}
 \int d\mv x\ \bar{\tau}_{(\mv x_1| \mv x_0)} (\mv x) &=& \int_0^\infty dT\ T\  p(\mv x_1 | \mv x_0) \\ 
 &=& p_\text{lin}(\mv x_1 | \mv x_0),
\end{eqnarray*}
that is, if we integrate the mean trail image, we just get the linearly reweighed version of \eqref{eq:linrew} and have
also a way to compute \eqref{eq:linrew} in practice.

\subsection{Application to DW-MRI data}
All concepts were introduced by considering a 3D diffusion process with drift. But, how to apply
this kind of stochastic process for 
brain connectivity analysis? The data is basically tensor-like, meaning it is point symmetric, there is no 
velocity field available. Even when we have orientation distributions, there is still no convection-like 
force. So, what to do? 

First of all, we have to realize that usual probabilistic tracking 
algorithms are not Markovian with respect to the position $\mv r\in\Real^3$. The propagation direction usually
depends not just on the position, but also on the previous step made. So, the state space 
of the tracker should be the joint space of position and orientation $\mv s = (\mv r,\mv n)$, where $\mv r\in \Real^3$
and $\mv n\in S_2$. The data we get from the diffusion measurement is basically an fiber orientation distribution 
$f(\mv r,\mv n)$ defined on this joint position/orientation space.
In fact, in this joint space there are several ways to define meaningful connectivity measures. 
 We propose to use the following Lagrangian
\begin{equation}
 L(\mv r,\td{\mv r},\mv n,\td{\mv n}) = \frac{1}{2\sigma_r^2} (\td{\mv r} - \mv n f(\mv r,\mv n))^2 + \frac{1}{2\sigma_n^2} \td{\mv n}^2 \label{eq:theLag}
\end{equation}
as path-costs. Paths with $\td{\mv s} \sim \mv n f$ and
small changes in $\mv n$ are assigned with small costs, and hence, with high probability.
 From the walker perspective,
there is a convective force which drives a walker with internal state $\mv n$ into direction $\mv n$, while the speed is proportional to the data term $f$. Additionally there is diffusion on
 the orientation variable, enabling the walker to change 
its directions, or conversely, penalizes too strong bending. 
Imagine the system as a network of pipes oriented in all possible directions carrying a fluid 
which flows with speed $f(\mv r,\mv n)$. The flow is not perfectly convective, 
neighboring parallel pipes can exchange fluid (this is spatial diffusion) and also 
crossing pipes exchange some fluid (this is the angular diffusion). 
In fact, we can also let spatial diffusion $\sigma_r \rightarrow 0$ such that 
there is only pure convection. However, as we will see later, the discretization of the FPE will force
us to have finite values for $\sigma_r$.

The FPE operator is not totally straight-forward, because there is an additional 
term coming from the curvature of the sphere (\cite{ito1981characterization,andersson1999finite}). It looks like
\[
 \mathcal{H} = \frac{1}{2} ( \sigma_r^2 \Laplace_{\mv r} + \sigma_n^2 \Laplace_{\mv n}) - \frac{1}{2} (\mv v \cdot \nabla_{\mv r} + \nabla_{\mv r} \cdot \mv v) + \frac{\sigma_n^2 R}{12}
 \]
where the velocity is defined by $\mv v(\mv r,\mv n) = \mv n f(\mv r,\mv n)$. The additional constant $\frac{R}{12}$ is 
coming from the curvedness of the sphere which is $R=2/\rho^2$ for a 2-sphere of radius $\rho$. However, there is actually 
no metric connection between the position space $\Real^3$ of the orientation space $S_2$. The 'radius' of the sphere $S_2$ has no actual meaning. 
So, the value of $R/12$ is rather arbitrary. In fact, $R/12$ acts like an additional exponential weighting, like the $\kappa$ we already introduced above
and could be absorbed by that. So, we exclude $R/12$ in the following, or, imagine $\rho$ to be pretty large. 

Let us now look at the symmetry properties 
of $\mathcal{H}$. Let $\mathcal{Z}$ the operator that does a point reflection on some function $a(\mv r,\mv n)$ in the sense
\[
 (\mathcal{Z} a)(\mv r,\mv n) := a(\mv r,-\mv n)
\]
Then, we can easily show that 
\[
 \mathcal{H Z} = \mathcal{Z H^+}
\]
because the data $f$ is symmetric $\mathcal{Z} f = f$. In other words, the operator
$\mathcal{H Z}$ is selfadjoint. Formally we can write the connectivity amplitudes as
\[
p(\mv r,\mv n|\mv r_0,\mv n_0) = -\langle e_{(\mv r,\mv n)}| \mathcal{H}^{-1} e_{(\mv r_0,\mv n_0)} \rangle
\]
where $e_{(\mv r',\mv n')}(\mv r,\mv n) = \delta(\mv r' -\mv r) \delta(\mv n' - \mv n)$.
For the connectivity amplitude the symmetry means
\[
p(\mv r,\mv n|\mv r_0,\mv n_0) = p(\mv r_0,-\mv n_0|\mv r,-\mv n),
\]
because $\mathcal{H}^{-1}  = \mathcal{Z}(\mathcal{H^+})^{-1}\mathcal{Z}$.

If we now think of some function $a(\mv r,\mv n)$ describing a seed point density,
and $b(\mv r,\mv n)$ a terminal point density. For both we can think of indicator
functions describing some cortical region of interests.
Then, we can define the symmetric connectivity measure $c(a,b)$ to be
\[
c(a,b) = -\langle a | \mathcal{Z H}^{-1} | b \rangle
\]
which can be written in terms of the connectivity amplitudes $p(\mv r,\mv n|\mv r_0, \mv n_0)$ 
as
\[
c(a,b) = \int \int a(\mv r,-\mv n) p(\mv r,\mv n|\mv r_0, \mv n_0) b(\mv r_0,\mv n_0)
\]
So, $c(a,b)$ is the path integral over all path connecting region $a$ with region $b$
according to the Lagrangian $L$. Accordingly, we can also define $c^\kappa$ and $c^\text{lin}$ corresponding
to  $p_\text{lin}$ in equation \eqref{eq:linrew}  and $p_\kappa$ in equation \eqref{eq:exprew}, respectively.
The path trail image is 
\begin{eqnarray*}
\bar{\tau}_{(a|b)}(\mv r,\mv n) 
&=& \langle a | \mathcal{H}^{-1} e_{\mv r,\mv n} \rangle \langle e_{\mv r,\mv n} | \mathcal{H}^{-1} b \rangle \\
&=&  \langle e_{\mv r,\mv n} | \mathcal{ZH}^{-1}\mathcal{Z} a \rangle \langle e_{\mv r,\mv n} | \mathcal{H}^{-1} b  \rangle \\
&=&  \langle e_{\mv r,-\mv n} | \mathcal{H}^{-1}\mathcal{Z} a \rangle \langle e_{\mv r,\mv n} | \mathcal{H}^{-1} b \rangle
\end{eqnarray*}
Note that $\mathcal{Z} \bar{\tau}_{(a|b)} = \bar{\tau}_{(\mathcal{Z}b|\mathcal{Z}a)}$.

Usually the seed functions $a$ and $b$ do not depend on $\mv n$
because we do not have any preference for the starting orientation, so $\mathcal{Z}a = a$ and
$\mathcal{Z}b = b$. Several relationships are then simplified, e.g. 
$\mathcal{Z} \bar{\tau}_{(a|b)} = \bar{\tau}_{(b|a)}$
and $c(a,b) = \langle a | \mathcal{H}^{-1} | b \rangle$.
For the special case if $a_{\mv r_0}(\mv r,\mv n) = \delta(\mv r-\mv r_0)$ we can also define the connectivity measure between two spatial locations $\mv r$ and $\mv r_0$ 
as 
\begin{equation}
c(\mv r_0,\mv r_1) = c(a_{\mv r_0},a_{\mv r_1})\label{eq:conampspat}
\end{equation}
In all experiments we will report the normalized connectivity amplitudes 
\[
 c_n(a,b) = \frac{c(a,b)}{\sqrt{c(a,a) c(b,b)}}
\]
which is quite natural and reminds one of Pearson's correlation coefficient. 


\subsection{Angular Constraints, the Data Term and Boundary Conditions}
The Lagrangian \eqref{eq:theLag} penalizes bending, which is a good prior for straight fibers, but underestimates curved 
fibers. In \cite{Behrens2007} this problem does not appear, because the walker is always back projected onto the nearest fiber direction in a voxel.  
Such a behavior cannot be realized by a Wiener process. However, we can mimic such a behavior 
by an angular constraint such that the walkers do not deviate too strong from the main fiber directions. 
That is, we only simulate paths where $\mv n$ is not too far away from the main fiber directions, which 
is similar to the maximum angle thresholds known from ordinary streamline algorithms, but not so rigorous (it
is still possible to find paths where the spatial tangent $\td{\mv r}$ to the path is not perfectly $\mv n$).
Now, we are able to set $\sigma_n$ pretty large, i.e. low penalty on bending, but the walker will be along the main fiber directions. 
In fact, this idea will also help to keep the running time and memory consumption of our algorithm in a reasonable range, because 
we only have to simulate the domain in the neighborhood of the fiber directions and not on the complete sphere. Note that the angular constraint
may also be seen as an additional cost in the Lagrangian, i.e. the forbidden path are assigned with infinite costs.

There are lots of ways to determine the speed function $f$, which represents the DW-data. For example, we could directly use 
the fiber orientation distributions (FOD), e.g. estimated by spherical deconvolution \cite{Tournier2007}. Or, one could estimate 
the main fiber directions by a local maxima detection of an FOD, and use them as anchor directions to construct the speed function. 
We opt for latter, primarily because the angular constraints can be much beter controlled. Let $\mv d_i$ be  the local
maximas of the FOD. Then, we construct the speed function by 
\begin{equation} 
f(\mv n) = \sum_{i=1}^N  (\mv n \cdot \mv d_i)^{2n} \label{eq:speedfun}
\end{equation} 
with some fixed $n>0$.
The simulation domain is now defined by thresholding this speed function, i.e. we only consider regions with $f>\epsilon$ for some
$\epsilon>0$. 

To solve the FPE we have to set some boundary conditions. It is natural to let the walkers die once they touch the boundary of the 
domain, which is equivalent to a zero Dirichlet conditions at the boundary. So the complete problem we want to solve is
\begin{equation} 
 -\mathcal{H} p = a \ \text{ where }\    p(\partial \Omega) = 0 \label{eq:problem}
\end{equation} 
where $\partial \Omega$ denotes the boundary of the simulation domain, $p$ the unknown steady state solution and $a$ the
seed region where the walkers are emitted. 
Of course, other boundary condition, like Neumann
conditions, might be applied, for example at the transitions to the ventricles. However, we restrict here all considerations
to the simple Dirichlet conditions.

\subsection{Implementation, Discretization and Solver}
In order to solve the equation $-\mathcal{H} a = b$
we have to discritize the operator $\mathcal{H}$. The convective part of $\mathcal{H}$ makes this a pretty hard problem, 
usual Finite Element Methods are only applicable for diffusion dominated problems. So, we decided to use a finite difference
upwind-downwind scheme. However, such schemes are known to introduce errors (known as false diffusion) when the 
convection direction is not aligned with the underlying discretization grid. In particular for steady state solutions these
problems become severe and the solutions show heavy orientation dependency. 
Fortunately our problem is special, the convection direction is always fixed for a volume $p(\mv r,\mv n)$ 
with fixed $\mv n$, namely $\mv n$. So, we propose to steer the grid for each direction along the current convection direction.
Practically, we discretize the sphere into $N$ direction homogeneously distributed over the full sphere (in the experiments we used
$N=128$ directions). So, we have $N$ volumes, where the discretization grid of the volume associated with the $i$th direction $\mv n_i$ is steered such that one of the coordinate axis is along direction $\mv n_i$. In this way we are able to produce quite clean convection, however the angular diffusion gets disturbed. A voxel at some discrete coordinate in volume $i$ does not have a unique partner in volume $j$ at exactly the same discrete position. We need this partner in order to implement the spherical Laplace-Beltrami operator $\Laplace_{\mv n}$, which only acts on the spherical coordinate. 
The simplest way to account for that is to use trilinear interpolation, which, unfortunately causes another unwanted spatial diffusion effect. However, this false diffusion is not so disturbing: it is independent of the orientation and the error 
is proportional to the angular diffusion and not to the speed function. Second, we can reduce the error by a finer discretization of the spatial grid. We quantified this additional spatial diffusion in an experiment and found that a angular displacement of $\pi = 180^\circ$  (for a $N=128$ points on the sphere) causes an additional spatial displacement of about $5$ spatial grid units, i.e. the cost of a $180^\circ$ turn on the spot compares to a 5 voxel jump. So, for some given theoretical $\sigma_n$ and $\sigma_r$, the actual 
$\sigma_r^\text{act}$ turns out to be approximately $\sigma^\text{act}_r = \sqrt{\sigma_r^2 + (\sigma_n 5/\pi)^2}$, if $\sigma_r$ 
is given in grid units. As this is already pretty much we decided to set the theoretical $\sigma_r$ to zero in all experiments. 
In general, for arbitrary number of points on the sphere, we found the approximation
$\sigma^\text{act}_r = \sqrt{\sigma_r^2 + (\sigma_n \sqrt N/7)^2}$, i.e. with a finer sphere discretization the 
effects get stronger. But again, note that the 'physical' amount of this false diffusion depends also on the spatial grid size. 
All the gritty details about the discretization scheme are given in Appendix \ref{app2}.

The whole matrix is scattered with MATLAB's sparse matrix capabilities. 
To solve the discretized equation \eqref{eq:problem} we found the GMRES algorithm the most stable and efficient \cite{saad1986gmres},
we just used the MATLAB's \verb+gmres+ command. As number of restarts we used $5$. The tolerance value is set to $10^{-6}$ in 
all experiments. The running time depends on the size of the problem. If we use the setting as described below in the 
experiments we get usually a system with about $2-3\cdot 10^6$ variables, which is scattered in $1-2$ minutes on a common Desktop PC (Intel I7, 16GB).
Solving the equation also takes about one minute.

\section{Experiments}
\begin{figure}
\begin{center}
\includegraphics[width = 9cm]{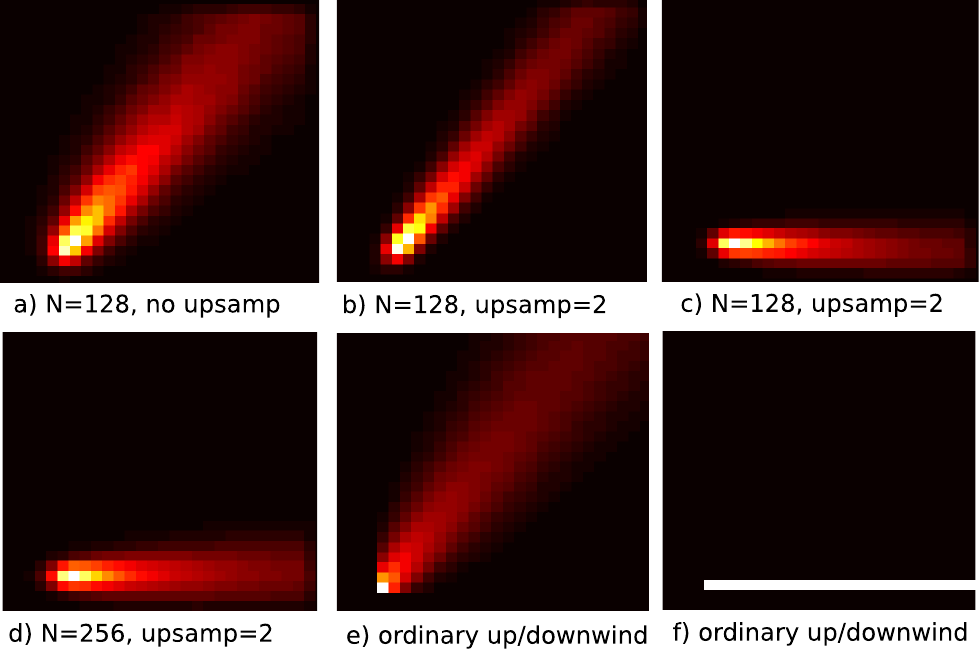} 
\end{center}
\caption{Demonstration of the proposed discretization scheme. In (a-d) the novel scheme for different granularities are shown. 
In (e-f) the ordinary scheme is shown for comparison. } \label{fig1}
\end{figure}

We want to begin with a small demonstration of the proposed discretization scheme. Let $\mv d$ some arbitrary 
direction and the speed function is constructed according to \eqref{eq:speedfun} with $n=25$ and $\epsilon = 0.02$. 
In Figure \ref{fig1} we choose $\mv d = (cos(\varphi),sin(\varphi),0)$ with $\varphi=0.3 \pi$ (a,b,e) and $\varphi=\pi$ (c,d,f). Further
we choose $\sigma_r=0$ and $\sigma_n=\pi/12$. Figure \ref{fig1} (a-d) shows the solution to \eqref{eq:problem} with $a(\mv r, \mv n) = \delta(\mv r-\mv r_0)$
where $\mv r_0$ is located in the lower left corner of the simulation domain. The spatial intensity maps in Figure \ref{fig1} show the solution $p$
integrated over the angular variable $\int_{S_2} p(\mv r,\mv n) d\mv n$, i.e. we actually show the connectivity amplitudes $c(\mv r,\mv r_0)$ as
defined in \eqref{eq:conampspat}.
The matrix size is of size $30\times30\times7$. 
In Figure \ref{fig1}a,b,c) we used a sphere with N=128 directions in (d) with N=256 directions. Further, we demonstrate the effect of spatial oversampling, i.e.
a upsampling by a factor of $2$ means that the actual matrix size is in the range of $60\times60\times14$ (for details see Appendix \ref{app2}). 
By comparing (a) and (b) we can obviously see how the the upsampling reduces false diffusion. 
Similarly comparing (b), (c) and (d), we see that the amount of false diffusion does not depend on the direction, but the directions are not perfect due to the discretization of the sphere.
For comparison we also simulated an ordinary up/downwind scheme in (e),(f) with $\sigma_n = \sigma_r = 0$. The amount 
of false diffusion depends dramatically on the direction.

\begin{figure}
\begin{center}
\includegraphics[width = 9cm]{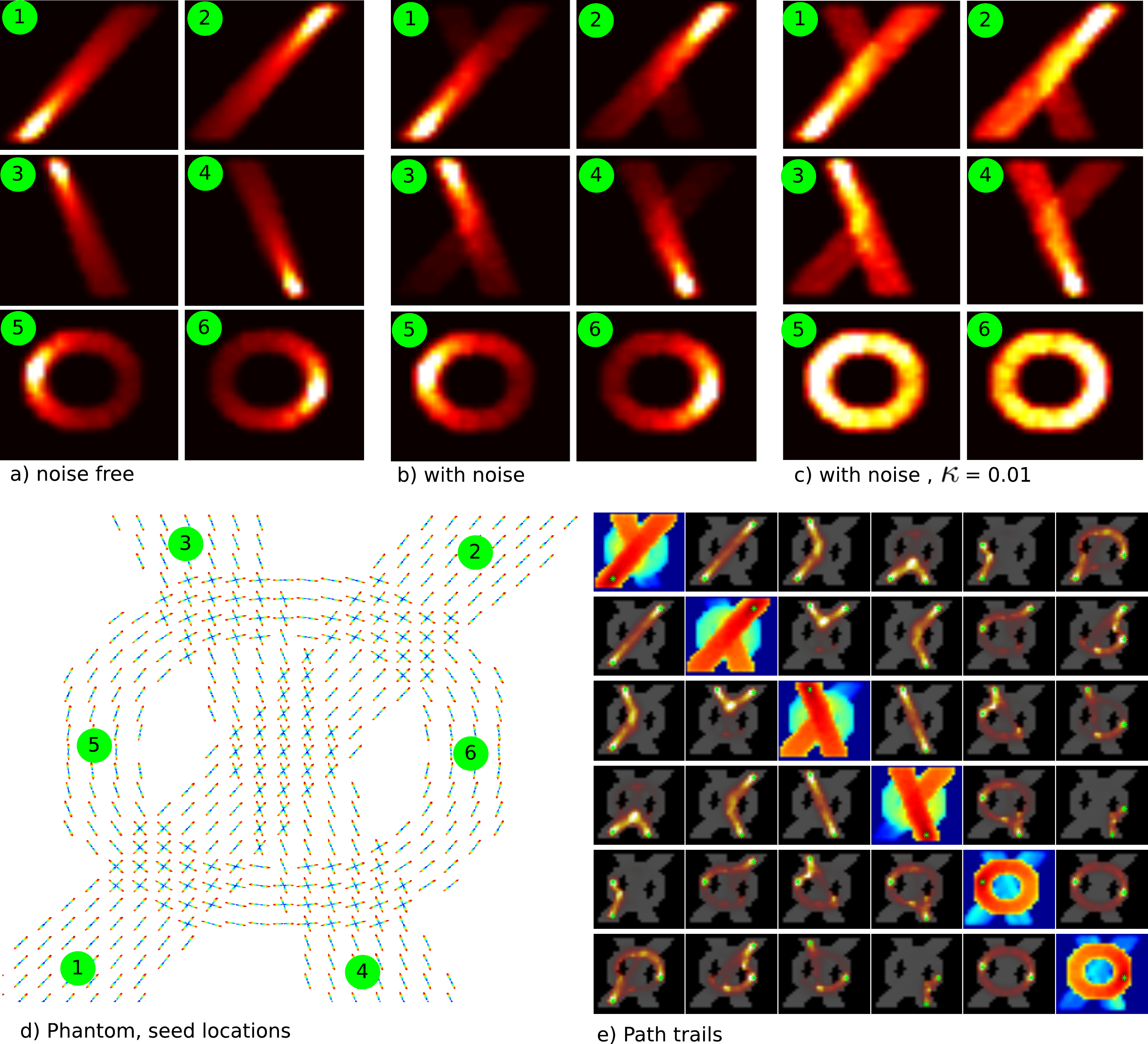} 
\end{center}
\caption{In (a-c) connectivity amplitudes for three different settings are shown, the color scaling is fixed and the same for all three settings.
In (d) the FOD of the phantom together with the seed locations is given, and (e) shows the matrix of all possible path trails together with 
the connectivity amplitudes in logarithmic scaling in the diagonal. 
} \label{fig2}
\end{figure}

\begin{figure*}
\begin{center}
\includegraphics[width = 18cm]{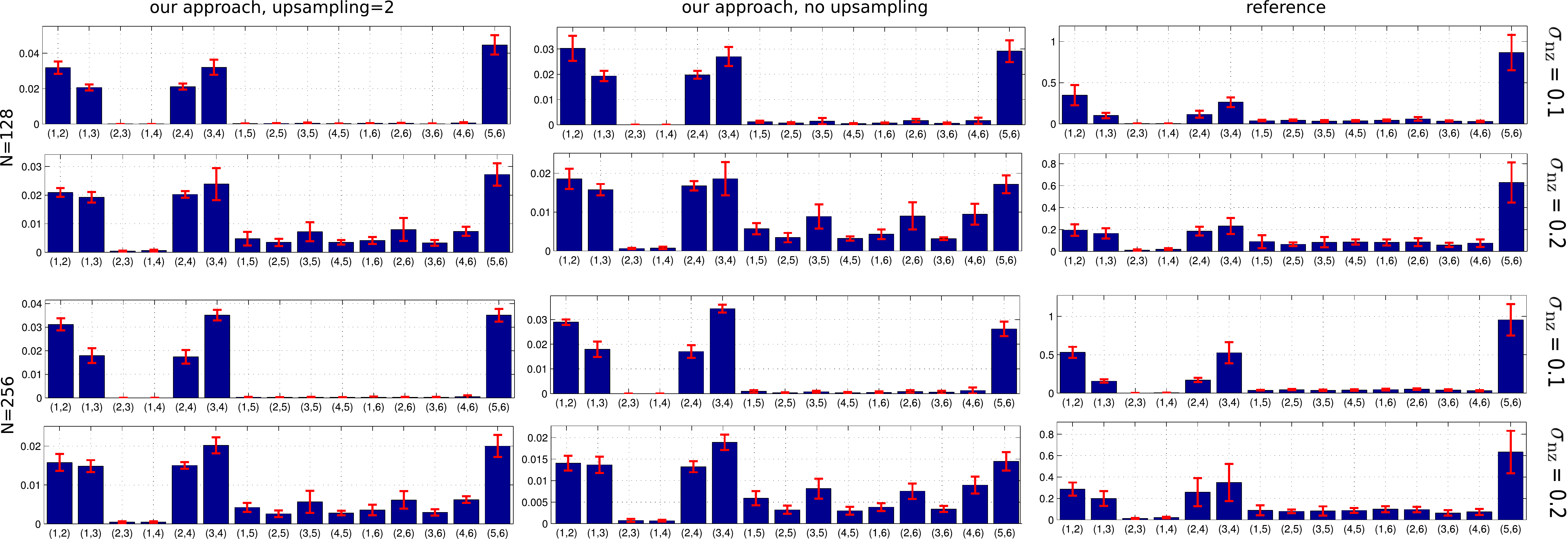} 
\end{center}
\caption{Quantitative analysis of connectivity values (c-values) over 50 trials. The barplots show mean c-values together with 
standard deviation. The x-axis indicate all possible seed pairs $(i,j)$ for $i,j=1,\hdots,6$ of the numerical phantom.} \label{fig3}
\end{figure*}

\subsection{Numerical Phantom}

For further demonstration we constructed a numerical phantom consisting of crossing and bending configurations (see Figure \ref{fig2}d)).
Six seed locations $\mv r_1,\hdots,\mv r_6$  were selected such that they are pairwise connected as (1-2), (3-4) and (5-6).
We did not simulate the MR-signal, but created the underlying directions directly.
The directions $\mv d_i$ are created continuously (not aligned with the discretized sphere directions), and a pseudo FOD is generated using 
$\sum_i \exp(\alpha ((\mv d_i \mv n)^2 - 1)$, which is shown in Figure \ref{fig2}d). 
For robustness analysis we distorted the directions by Gaussian noise of standard deviation $\sigma_{\text{nz}}$.

Figure \ref{fig2}a-c) gives some first results of our approach. We show the connectivity amplitudes 
$c(\mv r,\mv r_i)$ as images, where $\mv r_i$ is one of the fixed seed point. 
Figure \ref{fig2}a) shows $c$ without noise and b) with noise ($\sigma_{\text{nz}}=0.1$).
Figure \ref{fig2}c) shows the same experiment but with exponential length bias correction of $\kappa = 0.01$. Examples for the 
path trails images are shown in Figure \ref{fig2}e) that correspond to the experiment in Figure \ref{fig2}b).
Trail images for all pairs of seed combinations are shown, the diagonal shows again the connectivity amplitudes  but in 
logarithmic scaling. 

To get a quantitative picture we repeated the above experiment for different noise levels $\sigma_{\text{nz}} = 0.1,0.2$ and different discretization scheme,
namely spherical discretization with $N=128$ and $N=256$ directions and with and w/o spatial upsampling (with a factor of $2$). Additionally, 
the orientation of the discrete sphere was randomly changed for each run. As reference approach we followed \cite{Behrens2007} (for details see Appendix \ref{app3}),
with $\Delta s = 1, \sigma=0.2$, without revisits and with length bias correction. 
Figure \ref{fig3} shows barplots of the normalized connectivity amplitudes $c_n(\mv r_i,\mv r_j)$ for all pairs of seeds
averaged over 50 runs together with the standard deviation.

\begin{figure*}
\begin{center}
\includegraphics[width = 18cm]{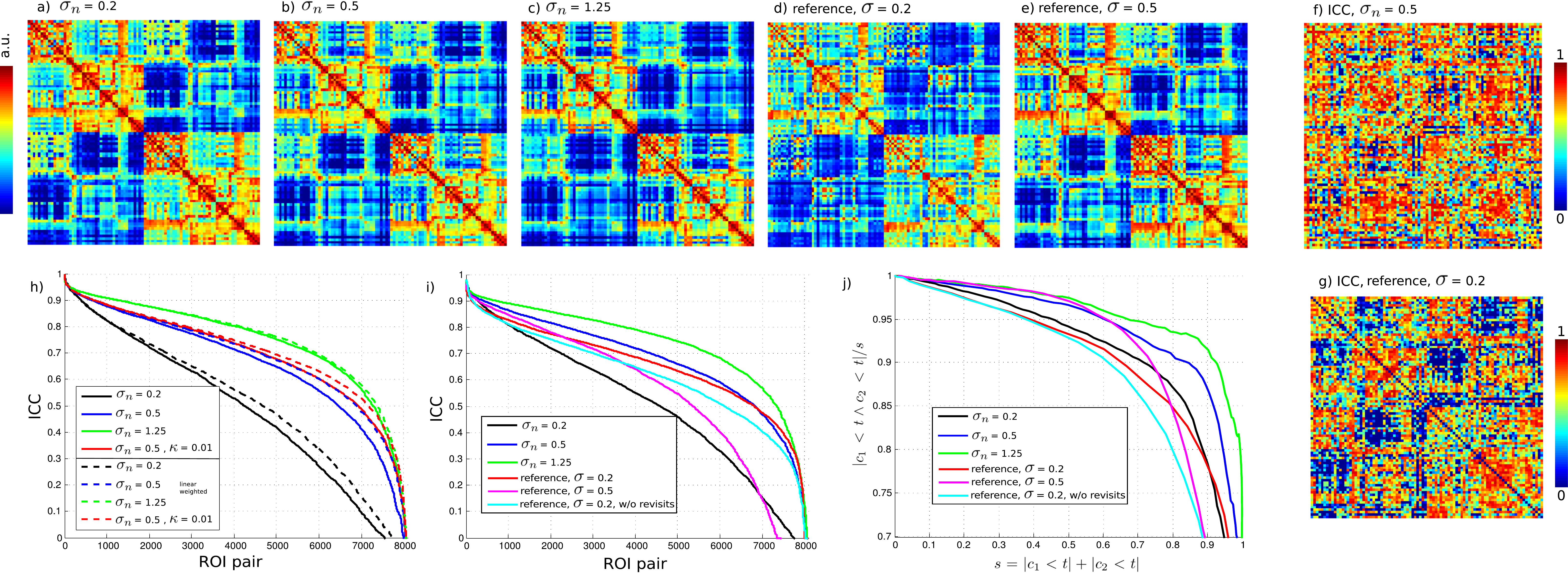} 
\end{center}
\caption{In (a-e) average connectivty matrices (CM) for different approaches are shown. The matrices in (f,g) show ICC values for individual 
c-values. Plot (h,i) display the ICC for all c-values (all ROI pairs) in the CM. Plot (j) shows agreement values.  } \label{fig4}
\end{figure*}

\subsection{In-Vivo Human Brain}
To investigate our approach on real DTI data we considered 28 scans of healthy volunteers at a b-value of $1 \text{ms}/\mu m^2$ with 61 diffusion directions and 
an isotropic resolution of $2mm^3$. A white matter probability map was generated with SPM (Version 8, \emph{http://www.fil.ion.ucl.ac.uk/spm/}, \cite{SPMbook2007})
on a T1-weighted scan, which was co-registered to the $b_0$-scan of the diffusion sequence. For each subject the scans were repeated two times (in two different
sessions) to allow to investigate the robustness of the approach. To compute the FOD a constrained spherical deconvolution was used together
with a $L_1$-regularizer \cite{Michailovich2011}. The fiber response function (FRF) was
chosen as $FRF(\mv n) = \exp(-b D_1 (n^\text{fib}_i n_i)^2)$ with $D_1 = 1 \mu m^2/ms$ for all subjects and reconstructed on a $N=128$ sphere.
The local maximas of the FOD are determined by fitting a quadratic form to the neighborhood of the discrete local maxima (the same 
approach as used in \cite{reisert2011}). Spurious local maximas were filtered out by neglecting all maximas smaller than $20\%$ of the global one. 
The speed function is constructed in the same way as for the numerical phantom. The spherical discretization is kept at $N=128$ and no spatial
up sampling is applied. The probabilistic white matter segmentation of SPM was thresholded at $0.5$ and used as white matter mask. 
To compare and analyze the group, the AAL atlas \cite{tzourio2002automated} was registered to the subjects' native space, normalization was done with SPM 8.
The first 90 ROIs (distributed over the neocortex) of the AAL atlas were used to compute connectivity matrices.

Figure \ref{fig4}a-e) show examples of the mean connectivity matrices (CM) over the whole group. The first 45 ROIs belong to the right 
hemisphere, the last 45 the right hemisphere. The CM is shown in logarithmic scaling. 
To get a comparable contrast we used the following formula $c_\text{log}(a,b) = \log(t+c_n(a,b))$ where $t$ is the $20\%$ quantile 
of $c_n$ over all regions. The CMs obtained by our approach are compared to the reference approach \cite{Behrens2007} (Figure \ref{fig4}d-e)).
To investigate the robustness we computed the intraclass correlation coefficient (ICC)\footnote{If $c_1$ and $c_2$ are c-values in the CM for scan 1
and scan 2, respectively, then the ICC is $ICC(c) = \langle(c_1-\bar c)(c_2-\bar c)\rangle / \langle(c-\bar c)^2\rangle)$, 
where $\langle\rangle$ denotes the expectation over the whole group and $\bar c = \langle (c_1 + c_2)/2 \rangle$. The ICC is 1 if all the
variance of the c-value can be explained by differences amongst the subjects.} for each connectivity value (c-value) in the CM. 
In Figure \ref{fig4}h) we show the ICC over all c-values in the CM sorted by ICC magnitude for different settings: we varied $\sigma_n$ and 
show one result with exponential weighting ($\kappa=0.01$). For all results we also show the corresponding linear weighted result $c_n^\text{lin}$. 
Figure \ref{fig4}i) compares our approach with \cite{Behrens2007} with different parameters: 
$\sigma = 0.2,0.5, \Delta s = 1$, with (w/o) revisits and length bias correction (lbc). 
In Figure \ref{fig4}f,g) we show for our approach and the reference the ICCs for the individual c-values as a matrix. 
Table \ref{tab1} shows mean and median ICCs for our approach and the reference with more different settings. 
It is common to put thresholds to get a binary decision of connectivity. To investigate the robustness of such a thresholding 
operation an agreement measure is calculated as follows: let $t$ be some threshold, then $a = |c_1<t \wedge c_2 <t|$ denotes the number of 
regions that are not connected in scan $1$ and in scan $2$, where the number is counted over all possible ROI pairs and subjects. The agreement value 
$a$ is normalized by the total number of non-connected regions $s = |c_1 <t| + |c_2<t|$, and hence, $a/s$ is a number between $0$ and $1$. If we 
now vary over all thresholds $t$ and plot $a/s$ as a function of $s$ we get plots displayed in Figure \ref{fig4}j). 

\begin{figure}
\begin{center}
\includegraphics[width = 9cm]{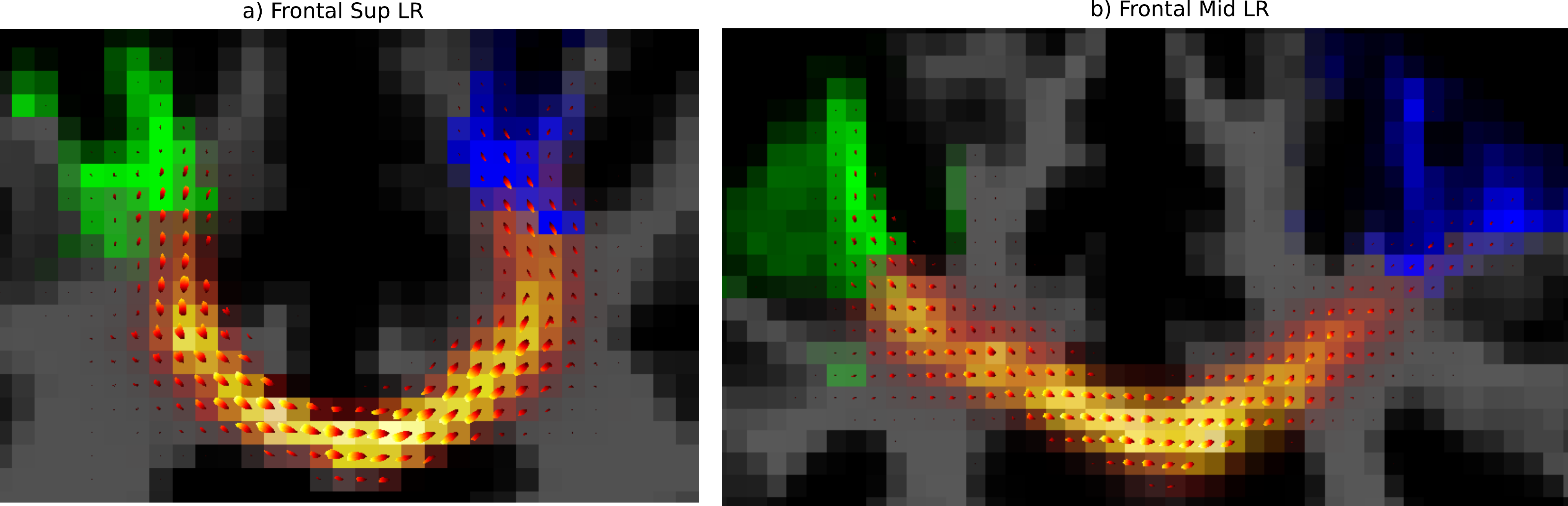} 
\end{center}
\caption{Examples of mean path trails for two ROI pairs in the frontal regions. Green and Blue depict the ROIs. 
The path trails are shown as glyphs and averaged over directions.} \label{fig5}
\end{figure}

Finally, Figure \ref{fig5} shows examples of path trail images for two ROI pairs of the AAL atlas. The two ROIs are highlighted in green and blue, respectively.
The white matter mask is depicted in dark gray. The path trail $ \bar{\tau}_{(a|b)}(\mv r,\mv n)$ is obviously
a function in the joint position/orientation space, which is shown in glyph representation, but note that the glyphs are not pointsymmtric as 
one used to observe. The traversal direction is in this example from the blue to the green ROI.

\begin{table}
\caption{Quantitative Comparison by ICC} \label{tab1}
\centering
\begin{tabular}{|l|c|c|}
\hline
\textbf{Method} & \textbf{Mean} & \textbf{Median} \\ 
\hline 
Our Approach \\
\hline
$\sigma_n = 0.2$ & 0.47  & 0.51  \\
$\sigma_n = 0.5$ & 0.64  & 0.69  \\                  
$\sigma_n = 1.25$ & 0.73  & 0.78  \\
$\sigma_n = 0.5, \kappa = 0.01$  & 0.67  & 0.72  \\
$\sigma_n = 0.2$, lin. weighted & 0.50  & 0.54  \\
$\sigma_n = 0.5$, lin. weighted  & 0.67  & 0.71  \\
$\sigma_n = 1.25$, lin. weighted & 0.73  & 0.78  \\
$\sigma_n = 0.5, \kappa = 0.01$, lin. weighted  & 0.69  & 0.73 \\
\hline 
Reference \cite{Behrens2007} \\
\hline
$\sigma=0.2, \Delta s = 1$  & 0.65 & 0.68 \\
$\sigma=0.5, \Delta s = 1$  & 0.55 & 0.63 \\
$\sigma=0.2, \Delta s = 1$,  w/o revisits  & 0.61 & 0.64 \\
$\sigma=0.2, \Delta s = 1$,  no lbc & 0.64 & 0.67 \\
$\sigma=0.2, \Delta s = 0.25$ & 0.55 & 0.58 \\
$\sigma=0.2, \Delta s = 0.5$ & 0.60 & 0.63 \\
\hline
\end{tabular}
\end{table}

\section{Discussion}
We proposed to use a Wiener process for the estimation of structural brain connectivities. The path integral perspective gave us the insight 
why the usual random walker principle leads to asymmetric connectivities. To symmetrize the solution an additional potential was introduced,
which symmetrized the path integral and the corresponding Fokker Planck equation. However, particle conservation and, hence, the usual 
walker perspective is lost. To solve the FPE we introduced a novel discretization 
scheme which avoids the directional dependency of the false diffusion introduced by the discretization of the convection operator. 

Initial experiments (in Figure \ref{fig1}) show that the new discretization approach is valid. In Figure \ref{fig2} the visual results on the numerical 
phantom show the approach works also nicely under noisy conditions and that the exponential length bias correction is reasonable.
Quantitative experiments show that it leads to a robust connectivity measure. Compared to \cite{Behrens2007}
the new approach shows smaller connectivity values for bending configurations, which is due to the bending penalty which is inherent to the new approach
(see Figure \ref{fig3}). One can also see that the robustness of our approach does not really depend too much on the discretization 
granularity. Compared to \cite{Behrens2007} it is in most cases higher (lower standard deviation over the trials).
For stronger noise ($\sigma_\text{nz}=0.2$) our approach gives relatively higher c-values to false connections, but the c-values itself are more stable.
Another difference to \cite{Behrens2007} is the way the curvature is controlled. If we want for our approach a higher connectivity for bending 
fibers (ROIs (5,6)), then we also get higher c-values for (1,3) and (2,4). For \cite{Behrens2007} it is different, the step directions 
are always back projected onto the nearest most likely fiber direction, which results in higher c-values for the bending (5,6).

To investigate the new approach in-vivo, we considered a group of volunteers which were all scanned twice. We did a ROI-based approach
based on the AAL-atlas. Overall, the CMs obtained from our approach and the reference are similar. 
As one might expect, with our approach strong bending fibers (like the connection between 'Occipital Sup L' to  'Occipital Sup R')
are assigned smaller connectivity values. 
Analysis of the intraclass correlations of c-values shows
that the new approach is quite robust. Usually higher angular diffusion $\sigma_n$ gives more robust c-values. Robustness is further increased
by a length bias correction, either linear or exponentially. Compared to \cite{Behrens2007} our approach shows good performance. For 
\cite{Behrens2007} increasing angular diffusion does not help so much. One can see that in regions with low connectivity \cite{Behrens2007} has problems
(see Figure \ref{fig4}f,g)), which goes inline with our finding in the numerical phantom. We found \cite{Behrens2007} to work best for $\Delta s = 1$, $\sigma = 0.2$,
revisits are allowed and with length bias correction (see Table \ref{tab1}). In this setting the median ICC is $0.68$ which is acceptable,
$50\%$ of the c-values have a ICC higher than $0.68$. The highest value we can achieve with our approach is a median of $0.78$. 
The agreement scores obtained from thresholded c-values show a similar but more pronounced picture (see Figure \ref{fig4}j)).
From Figure \ref{fig5} it can also be seen that the estimated path trail images are meaningful and follow the expected pattern.

%
%
%

\section*{Appendix}

\subsection{Operator Formalism} \label{app1}
The space of square integrable function on $\Real^n$ is a Hilbert 
space $H$ with the inner-product
\[
\langle q|p \rangle = \int_{\Real^n} \conj{q}(\mv x) p(\mv x) d\mv x
\]
Sometimes we adopt Dirac's bra-ket notation, i.e. by $|p\rangle$ we 
denote a vector in $H$, by $\langle p|$ a linear 1-form.
The adjoint of a linear operator $\mathcal{A}$ is
defined to be the operator $\mathcal{A}^+$ such that 
$\langle q|\mathcal{A} p \rangle = \langle \mathcal{A^+} q|p \rangle$
 for all functions $q,p \in H$. 

\subsection{Discretization} \label{app2}
Let $m(\mv r)$ some predefined white matter mask. For each direction $\mv n_i$ on the sphere, we create a voxel grid $\mv r^i_{(x,y,z)}$ such that 
the $x$-axis is aligned with direction $\mv n_i$. Only valid triples $(x,y,z) \in \Natural^3$ are considered where $m(\mv r^i_{(x,y,z)})$ is on. 
The voxel spacing $h = |\mv r^i_{(x+1,y,z)}-\mv r^i_{(x,y,z)}|$ can be chosen arbitrarily. 
In the experiments we used $h=2mm$ which is resolution of the DTI scan, and $h=1mm$, that is an oversampling factor of 2.
The sphere has a natural neighborhood system $\mathcal{N}(i)$
obtained by Voronoi tessellation. Let $p_{i,x,y,z}$ be some function on the introduced discrete domain. 
To discretize the spherical Laplace-Beltrami $\Laplace_n$ operator we define the discrete operator $\tilde \Laplace_n$ as
\[
 (\tilde \Laplace_n p)_{i,x,y,z} = \sum_{k\in \mathcal{N}(i)} \left( p_{i,x,y,z} -\sum_{(a,b,c)\in \atop \mathcal{T}(\mv r^k_{(x,y,z)})} w_{k,a,b,c}\  p_{k,a,b,c}  \right)
\]
where $\mathcal{T}(\mv r^k_{(x,y,z)})$ is neighborhood system needed for trilinear interpolation at position $\mv r^k_{(x,y,z)}$ in the $k$-th volume,
and $w_{k,a,b,c}$ are the corresponding weights. However, the matrix $\tilde \Laplace_n$ is not symmetric in the sense of the matrix transpose due to 
the one-sided interpolation. So the final discretized spherical Laplacian  is $\frac{\sigma_n^2}{2A}(\tilde \Laplace_n + \tilde \Laplace_n^\T)$. The
factor $A$ is empirically determined to get standard deviations of approximately $\sigma_n$.

To implement the convection part, the speed function $f(\mv r,\mv n)$ given in the native DTI frame is also transferred by  trilinear interpolation onto the new coordinate 
frame $f_{i,x,y,z} = f(\mv r^i_{(x,y,z)},\mv n_i)$. The simulation domain is determined by thresholding $f_{i,x,y,z}$. For the convection operator 
we have to distinguish between 'real' boundaries and boundaries coming from the thresholding.
The discretized convection generator $\tilde \nabla$ is defined as
\[
 -(\tilde \nabla p)_{i,x,y,z} = \left\{ \begin{array}{ll} \frac{p_{i,x-1,y,z} - p_{i,x,y,z}}{h} & \text{otherwise}\\
                                 p_{i,x-1,y,z}/h &  \text{if }b_{i,x+1,y,z}=1 \\
                                 -p_{i,x,y,z}/h &  \text{if }b_{i,x-1,y,z}=-1 \\
                                 -p_{i,x,y,z}/h &  \text{if }b_{i,x-1,y,z}=1 \\ 
                                \end{array} \right.
\]
where $b_{i,x,y,z}$ is an boundary indicator, where $b=1$ means the boundary is coming from thresholding and $b=-1$ from a 'real' boundary.
Finally, the convection part of the equation is $-(F \tilde \nabla + \tilde \nabla F)$, where $F$ is the matrix with $f_{i,x,y,z}$ on the diagonal.

\subsection{Reference Method} \label{app3}
For comparison we followed \cite{Behrens2007}. Suppose for some given fiber directions, which are always extracted in the same way as for our
approach. The state of a walker is its current position $\mv r$ and its last step direction $\mv n$. The 
propagation law for a single walker is: The current voxel contains a set of fiber directions $\mv d_i$. Determine the fiber direction $\mv d$ which 
is closest to the current $\mv n$. If the angle between those is above a threshold of $80^\circ$ tracking is stopped.
Otherwise, disturb this fiber direction by Gaussian noise $\mv d_n = \mv d + \eta_\sigma$ with deviation $\sigma$. Renormalize $\mv d_n$
to unity and track along this direction with step width $\Delta s$, i.e. $\mv r^\text{new} = \mv r + \Delta s \mv d_n$ and set $\mv n^\text{new} = \mv d_n$. If the new
position is out of the tracking mask, stop tracking. If desired, one can also stop tracking whenever a voxel is visited more than once. During 
tracking we count how often a walker has visited a voxel resulting in a probabilistic map (PM). Alternatively, to account for the length bias,
one can also weight with visitation count by the mean length of all walkers visited the voxel. To compute the connectivity value between ROI A and 
ROI B, $N$ walkers are ejected in each voxel in ROI A.
The connectivity value is then the sum over the PM in ROI B. To get a symmetric value the 
connectivity from A to B is averaged with the connectivity from B to A. 
To set $N$ we computed the intra-scan ICC, i.e. for each subject \emph{and} scan 
in the group we started the algorithm twice and computed the ICC between the two algorithm runs. 
We set $N$, such that the median ICC over all ROI pairs is above $0.95$, which resulted in $N=5000$ per voxel.

 \bibliographystyle{plain}

\end{document}